\documentclass[final,3p,times]{elsarticle}

 \usepackage{graphics}
 \usepackage{graphicx}
 \usepackage{epsfig}
\usepackage{caption}
\usepackage{subcaption}
\usepackage{varwidth}
\usepackage{multirow}
\usepackage{hhline}
\usepackage{wrapfig}
\usepackage[boxruled]{algorithm2e}

\usepackage{ulem}
\usepackage{algorithmic}
\usepackage{graphicx}
\usepackage{geometry}

\DeclareCaptionFormat{myformat}{\fontsize{10}{15}\selectfont#1#2#3}
\usepackage{color} 
\usepackage{amsmath}
\usepackage{amssymb}
\usepackage{amsthm}

\newcommand{\s}{{\mathbf s}}
\newcommand{\nn}{{\mathbf n}}

\def\mb{\mathbf}

\journal{NA}

\begin{document}

\begin{frontmatter}



\title{ Accurate adaptive deep learning method for solving elliptic problems}
\author{Jinyong Ying\fnref{label1}}
\author{Yaqi Xie\fnref{label1}}
\author{Jiao Li\fnref{label2}}
\author{Hongqiao Wang\fnref{label1}\corref{cor1}}
\fntext[label1]{School of Mathematics and Statistics, HNP-LAMA, Central South University, Changsha, Hunan 410083 China}
\fntext[label2]{School of Mathematics and Statistics, Hunan Provincial Key Laboratory of Mathematical Modeling and Analysis in Engineering, Changsha University of Science and Technology, Changsha, Hunan 410114, China}
 \cortext[cor1]{Corresponding author, email address: Hongqiao.Wang@csu.edu.cn.}
\begin{abstract}
 Deep learning method is of great importance in solving partial differential equations. In this paper, inspired by the failure-informed idea proposed by Gao  \textit{et.al.} (SIAM Journal on Scientific Computing 45(4)(2023)) and as an improvement, a new accurate adaptive deep learning method is proposed for solving elliptic problems, including the interface problems and the convection-dominated problems. Based on the failure probability framework, the piece-wise uniform distribution is used to approximate the optimal proposal distribution and an kernel-based method is proposed for efficient sampling.
Together with the improved Levenberg-Marquardt optimization method, the proposed adaptive deep learning method shows great potential in improving solution accuracy. Numerical tests on the elliptic problems without interface conditions, on the elliptic interface problem, and on the convection-dominated problems demonstrate the effectiveness of the proposed method, as it reduces the relative errors by a factor varying from $10^2$ to $10^4$ for different cases.
\end{abstract}


\end{frontmatter}


\section{Introduction}
Solving partial differential equations (PDEs) is important in applications and simulations. Due to powerful approximation abilities of neural networks, deep learning methods gain mathematicians' attention. The method was proposed and had been used since last century \cite{lagaris1998artificial,lagaris2000neural}, and became more and more polular as it being able to solve the Poisson equation using the deep Ritz method \cite{weinan2018deep} and other PDEs (as well as other related inverse problems) using the typical physics-informed neural networks (PINNs) \cite{raissi2019physics}. These two methods are similar to other other, just using different loss functions. More specifically, constructions of loss functions for the deep Ritz method are based on weak formulations of symmetric elliptic problems, while loss functions for the PINNs are based on the residual forms of the original PDEs. Due to the requirement of symmetry in construction of the loss function, PINNs are more widely used for solving PDEs, compared to the deep Ritz method. As for elliptic interface problems, these two methods encountered low accuracy issue. To overcome this difficulty, there are mainly two types of improvements. One type of improvements is 
using the domain decomposition idea. That is, to approximate the solution of the interface problem, two or more neural networks instead of a single network are used to simulate the solution within each subdomain. Here multiple networks are `connected' by enforcing the interface conditions from different subdomains. Using either the residual-based loss function \cite{he2022mesh} or the weak type loss function \cite{guo2022deep}, the improved deep learning method is shown to be able to get approximations with improved solution accuracy (relative errors being about $10^{-3}$) using the Adam solver. The other type is based on the extension-and-projection technique. That is, by augmenting one more dimensionality in spatial variable \cite{hu2021discontinuity}, one network is used to approximate the solution of the elliptic interface problem. Meanwhile, numerical experiments show that using the Levenberg-Marquardt (LM) method instead of the Adam method, the numerical solution can reach quite high accuracy (relative errors being about $10^{-6}$ or even better for some cases). Further study \cite{hou2023hrw} shows that the improved deep learning method based on the domain decomposition idea can also reach high solution accuracy, although using one network can achieve slightly better results than multiple networks. 
Besides, for the deep learning method, many improvements with respect to different aspects have been proposed in order to gain more accurate approximations, including improving the topology structure of the neural network \cite{paszke2016enet,ying2022multi} to enhance its representation ability, embedding the gradient information of PDEs \cite{yu2022gradient}, and optimally choosing the weights in loss functions \cite{wu2022inn}.

As for solving elliptic problems with complex solution structures (e.g., solution having a `peak' point or containing the internal/boundary layer), how to effectively choose the training points is significant since a fixed set of (uniformly chosen) training points may fail to capture the effective solution region \cite{wu2023comprehensive,wang2022and}, which leads to the proposals of adaptive deep learning methods based on different strategies \cite{gao2023failure,lu2021deepxde,tang2021deep}. Among these, the residual-based adaptive refinement method \cite{lu2021deepxde} is the most commonly-used strategy, which adds the sampling points simply based on predictions of residuals of  loss functions. Another interesting one is the so-called deep adaptive sampling (DAS) method \cite{tang2021deep}, where a generative model is used to generate new training points based on the residuals. The last one is just the failure-informed PINNs (FI-PINNs) \cite{gao2023failure}, which novelly introduces a failure probability estimation and uses a truncated Gaussian-based approximation  to automatically generate new training points. More specifically, the FI-PINNs employ an importance sampling strategy to evaluate failure probabilities. This method can not only allocate sampling points selectively for better training of PINNs, but also intelligently determine the termination of training point increments. Meanwhile, the importance sampling technique is anticipated to improve efficiency compared to the Monte Carlo (MC) method. However, it necessitates the density expression of the optimal proposal distribution. FI-PINNs address this requirement by approximating the optimal proposal distribution using either a truncated Gaussian distribution or a mixture Gaussian distribution with an adaptive scheme, enhancing the robustness of the methodology. As shown in \cite{gao2023failure}, the FI-PINNs is far more effective than the residual-based adaptive refinement method in terms of solution accuracy and is better than the DAS method in terms of efficiency. 

In this paper, as an improvement, we are trying to propose a new method which employs the MC method to estimate the failure probability. While this approach is computationally less efficient, it offers the advantage of bypassing the need of the explicit expression of the optimal proposal estimation. Meanwhile, the Monte Carlo simulations are cost-effective due to the utilization of a neural network architecture, and the enhanced sampling efficiency  from the optimal proposal can be achieved by further using the approximation of a piece-wise uniform distribution. 

Based on these techniques, we then carry out numerical experiments, where several different types of elliptic equations are considered. First of all, the problem with one peak in solution is tested and used to show the improvements of the new adaptive sampling method in terms of the self-adaptive importance sampling (SAIS) method used in FI-PINNs. In section \ref{example1}, the numerical results demonstrate that using the new method, the relative $L^2$ error not only reduces to $8.41\times 10^{-6}$ from $1.44\times 10^{-1}$, better than the one ($4.98\times 10^{-4}$) obtained using the SAIS method, but also decreases much faster. Similar results are obtained for problem with two peaks in solution and the high dimensional test problem. Then, considering one elliptic interface problem and two convection-dominated problems, the numerical results also demonstrate that the method still works by further reducing the relative errors. All numerical results show the effectiveness of the proposed adaptive sampling method.

The paper is organized as follows. Firstly, considering the elliptic problem with the interface condition or not,  the deep learning methods are introduced in Section 2, followed by the optimization method as well as  the proposed adaptive sampling method. Then various examples are used to show the effectiveness of the new method in Section 3. And finally a conclusion is given.

\section{Deep learning method}
In this section, we consider the following general elliptic PDE to introduce the details of the deep learning method. Consider the following general elliptic PDEs:
\begin{equation}
\left\{\begin{array}{ll}
-\nabla\cdot (\alpha(x)\nabla u(x)) +\beta(x)\cdot\nabla u(x) = f(x),  &\quad x \in \Omega, \\
    u(\s^+)=u(\s^-)+p(\s), \quad
 \displaystyle \alpha(\s^+)\frac{\partial u(\s^+)}{\partial\nn(\s)}=\alpha(\s^-)\frac{\partial u(\s^-)}{\partial\nn(\s)}+q(\s), &\quad  \s\in\Gamma,\\
u(\s) = g(\s), &\quad \s\in\partial\Omega,
\end{array}\right.
\label{poisson}
\end{equation}
where $\alpha(x)$ denotes the diffusion coefficient, $\beta(x)$ is a given vector function, $\Omega=\Omega_p\cup\Omega_s\cup\Gamma$ with the possible interface $\Gamma$ being the intersection of the subdomains $\overline{\Omega}_p$ and $\overline{\Omega}_s$ is a bounded convex domain in $\mathbb{R}^d$ ($d\geq 1$ is a specified positive integer), $g(x), p(\s)$ and $q(\s)$ are three given functions defined on the domain boundary $\partial\Omega$ and on the interface $\Gamma$, respectively, and $\nn$ is the outward unit normal vector of the domain $\Omega_p$.

For the general elliptic problem \eqref{poisson}, it contains several different types of elliptic problems considered in the literature. When the coefficient $\alpha(x)$ is continuously differentiable and $p(x)=q(x)=0$, it is just a diffusion or convection-diffusion problem without any interface, depending on magnitudes of the vector function $\beta(x)$. When the coefficient $\alpha(x)$ is assumed to be the simplest piece-wise constant and has jumps across the interface $\Gamma$, i.e.,
\begin{equation}
\alpha(x) = 
\left\{\begin{array}{cl}
\alpha_1,& {\text { if } x \in \Omega_p,}  \\
\alpha_2,& {\text { if } x \in \Omega_s,}
\end{array}\right.
\end{equation}
then the equation \eqref{poisson} becomes an elliptic interface problem. For these different types of elliptic problems, using the classic numerical methods such as the finite element method, formulations of weak forms and thus the discretization schemes are totally different in order to get accurate numerical solutions. Here we try to design a unified deep learning method to solve different types of problems and hope to get satisfactory results. Using the deep learning method for solving the elliptic problem with or without interface,  in order to get satisfying numerical solutions, the setups are slightly different. Hence in the following we divide into two cases to consider.

\subsection{Elliptic problems without interface conditions}
When we consider the  elliptic problem without any interface condition, for simplicity,  the problem can be rewritten as follows:
\begin{equation}
\left\{\begin{array}{ll}
-\nabla\cdot (\alpha(x)\nabla u(x)) +\beta(x)\cdot\nabla u(x) = f(x),  &\quad x \in \Omega, \\
u(\s) = g(\s), &\quad \s\in\partial\Omega.
\end{array}\right.
\label{poissonwithoutinterface}
\end{equation}

Here we use the deep learning method to solve the equation \eqref{poissonwithoutinterface} and try to set up the loss function, a key ingredient in the deep learning method. For simplicity, 
here the loss function is set to be the following residual-based or the least square loss function to characterize the solution $u$ of the elliptic problem \eqref{poissonwithoutinterface}:
\begin{equation}
\label{cminnointerface}
\mathop{\min}_{\tilde{u}\in H^2(\Omega)} L(\tilde{u}) =  \mathop{\min}_{\tilde{u}\in H^2(\Omega)}\sum_{i=1}^2\beta_i L_i(\tilde{u}),
\end{equation}
where $\beta_i$, $i=1,2$, are some given positive weight constants, and $L_i(\tilde{u})$, $i=1,2$, are the functionals defined as follows:
\[L_1(\tilde{u}) = \int_{\Omega}\left|-\nabla\cdot (\alpha(x)\nabla \tilde{u}(x)) + \beta(x)\cdot \nabla \tilde{u}(x)- f(x)\right|^2, \quad L_2(\tilde{u}) = \int_{\partial\Omega}|\tilde{u}(\s) - g(\s)|^2.\]

To minimize the loss function,  we firstly need to evaluate the integrals in the loss function, which can be effectively and efficiently approximated by the Monte Carlo (MC) algorithm by summing random
points in the corresponding domain. Introducing the uniformly distributed random points or the Latin hypercube sampling points \cite{loh1996latin}, $\left\{ x_i \right\}_{i = 1}^{{M}} \in {\Omega },\left\{ x_i \right\}_{i = 1}^{{M_{\partial \Omega }}} \in \partial \Omega $, in each subdomain, the corresponding optimization task becomes:
\begin{equation}
\label{dlf_noi}
{\tilde u^*_d} = \mathop {\arg \min }\limits_{\tilde u \in {\mathcal{H}}} {L }\left( {\tilde u} \right) = \mathop {\arg \min }\limits_{\tilde u \in {\mathcal{H}}} \sum\limits_{i = 1}^2 {{\beta _i}{l_i}\left( {\tilde u} \right)}, 
\end{equation}
where the hypothesis space $\mathcal{H}$ denotes the set of all expressible functions by a network with fixed width and depth given in the following
\begin{equation}
	\mathcal{H} = \left\{ {\left. {\left( {\cdot ,\theta } \right):{\mathbb{R}^{{d_{in}}}} \to {\mathbb{R}^{{d_{out}}}}} \right|\theta  \in \Theta } \right\} \nonumber
\end{equation}
with $\Theta$ denoting the undetermined parameters in the neural network, and the functionals ${l_i}\left( {i = 1,2} \right)$ are the discrete forms of $L_i$ given by
\begin{equation}
	{l_1}\left( {\tilde u} \right) = \frac{{\left| {{\Omega }} \right|}}{{{M}}}\sum\limits_{i = 1}^{{M}} {{{\left| { - \nabla \left( \alpha (x_i) \nabla \tilde u( x_i ) \right) + \beta(x_i)\cdot\nabla \tilde u\left(x_i \right)- f\left( {{{x}_i}} \right)} \right|}^2}}, \quad
	{l_2}\left( {\tilde u} \right) = \frac{{\left| {\partial \Omega } \right|}}{{{M_{\partial \Omega }}}}\sum\limits_{i = 1}^{{M_{\partial \Omega }}} { {{{\left| { {\tilde u\left( {{{s}_i},\phi \left( {{{s}_i}} \right)} \right) - g\left( {{{s}_i}} \right)} } \right|}^2}} } \nonumber
\end{equation}
with $\left| \cdot \right|$ being the measure of the corresponding subdomain.

Based on the simplest full feedforward neural network or the residual neural network, this method can produce numerical solutions with acceptable accuracy when the magnitude of the vector function $\beta$ is small (i.e., the diffusion problem) and the right-hand side function $f(x)$ is properly taken. However, when the convection dominance, the `peak' point or the internal/boundary layer is presented, the accuracy is limited. Combined with other improvements, we will overcome this drawback in the following.

\subsection{Elliptic problems with interface conditions}
For elliptic problems with interface conditions, directly applying the deep learning methods described before can not get accurate numerical solution. To overcome this accuracy issue, several techniques were proposed, such as the domain decomposition method \cite{he2022mesh,guo2022deep,ying2022multi} and the extension-and-projection method \cite{hu2021discontinuity, lai2022shallow, tseng2023cusp, li2023continuity}. Here instead of using the domain decomposition method,  we adopt the extension-and-projection method, since our recent work \cite{hou2023hrw} found that using one network based on this extension-and-projection method can achieve slightly better results than using multiple networks, and meanwhile it can preserve the continuity of the underlying quantity if necessary.

Given a neural network, to guarantee the accuracy, according to the extension-and-projection method \cite{hu2021discontinuity}, we firstly make the extension to one-higher spatial dimension augmented with the indicator function $z(x)$. That is, we try to approximate $\widetilde u(x,z)$ defined on $\widetilde \Omega   \equiv \Omega  \times \mathbb{R} \in {\mathbb{R}^{d + 1}}$ satisfying
\[ u(x) = \widetilde u\left(x,  z(x)\right),\]
where the indicator function $z(x)$ is defined to distinguish different subdomains $\Omega_p$ and $\Omega_s$. For the discontinuous problem on the interface, we can take the indicator function $z(x)$ to be the following piecewise constant
\begin{equation}
z(x) = 
\left\{\begin{array}{cl}
-1,& {\text { if } x \in \Omega_p,}  \\
1,& {\text { if } x \in \Omega_s.}
\end{array}\right.
\end{equation}
For the interface problem having continuous function values across the interface $\Gamma$, the indicator function $z(x)$ can be taken to $z(x)=\left| \psi(x) \right|$, where $\psi(x)$ is the level set function characterizing the interface of the underlying problem. As shown in \cite{li2023continuity}, this extension-and-projection deep learning method not only can preserve the continuity of the numerical solution on the interface $\Gamma$, but also can slightly improve the accuracy. Due to these advantages, we here use this `improved' method when the interface problem has continuous solution across the interface. 

Based on the extension function, we then try to construct the following loss function to characterize the solution $u$ of the elliptic interface problem \eqref{poisson}:
\begin{equation}
\label{cmin}
\mathop{\min}_{\widetilde{u}\in H^{3}(\widetilde\Omega)} L(\widetilde{u}) =  \mathop{\min}_{\widetilde{u}\in H^{3}(\widetilde\Omega)}\sum_{i=1}^5\beta_i L_i(\widetilde{u}),
\end{equation}
where $\beta_i$, $i=1,\cdots,5$, are some given positive weight constants, and $L_i(\tilde{u})$, $i=1,\cdots,5$, are the functionals defined as follows:
\[L_1(\widetilde{u}) = \int_{\Omega_p}\left|-\nabla\cdot (\alpha(x)\nabla \widetilde u\left(x,  z(x)\right)) + \beta(x)\cdot \nabla \widetilde u\left(x,  z(x)\right)- f(x)\right|^2, \]
\[L_2(\widetilde{u}) = \int_{\Omega_s}\left|-\nabla\cdot (\alpha(x)\nabla \widetilde u\left(x,  z(x)\right)) + \beta(x)\cdot \nabla \widetilde u\left(x,  z(x)\right)- f(x)\right|^2, \]
\[L_3(\widetilde{u}) = \omega_p \int_{\Gamma}\left|  \widetilde u\left(\s^+,  z(\s^+)\right) - \widetilde u\left(\s^-,  z(\s^-)\right) - p(\s) \right |^2, \quad L_4(\widetilde{u}) = \int_{\Gamma}\left|  \alpha(\s^+)\frac{\partial u\left(\s^+,  z(\s^+)\right)}{\partial\nn(\s)} - \alpha(\s^-)\frac{\partial u\left(\s^-,  z(\s^-)\right) }{\partial\nn(\s)} - q(\s)   \right|^2, \]
\[L_5(\widetilde{u}) = \int_{\partial\Omega}| \widetilde u\left(\s,  z(\s)\right) - g(\s)|^2.\]
Here $\omega_p$ is a determined parameter such that $\omega_p=0$ when $p(\s)\equiv 0$ and $\omega_p=1$ otherwise, used to distinct the considered two cases.

Same as before, the integrals in the loss function are effectively and efficiently approximated by the MC algorithm by summing random
points in the corresponding subdomain. Introducing the uniformly distributed random points, $\left\{ x_i \right\}_{i = 1}^{{M_p}} \in {\Omega_p}, \left\{ x_i \right\}_{i = 1}^{{M_s}} \in {\Omega_s }, \left\{ x_i \right\}_{i = 1}^{{M_{\Gamma}}} \in {\Gamma }, \left\{ x_i \right\}_{i = 1}^{{M_{\partial \Omega }}} \in \partial \Omega $, in each subdomain, the discrete optimization problem becomes:
\begin{equation}
\label{dlf_i}
\widetilde u^*_d=\mathop {\arg \min }\limits_{\widetilde u \in {\mathcal{H}}} {L }\left( {\tilde u} \right) = \mathop {\arg \min }\limits_{\widetilde u \in {\mathcal{H}}} \sum\limits_{i = 1}^5 {{\beta _i}{l_i}\left( {\widetilde u} \right)}, 
\end{equation}
where $\mathcal{H}$ again denotes the set of all expressible functions given by
\begin{equation}
	\mathcal{H} = \left\{ {\left. {\left( {\cdot ,\theta } \right):{\mathbb{R}^{{d_{in}+1}}} \to {\mathbb{R}^{{d_{out}}}}} \right|\theta  \in \Theta } \right\} \nonumber
\end{equation}
with $\Theta$ denoting the undetermined parameters, and the functionals ${l_i}$ are the discrete forms of $L_i$ given by
\begin{align*}
	{l_1}\left( {\widetilde u} \right)& = \frac{{\left| {{\Omega_p }} \right|}}{{{M_p}}}\sum\limits_{i = 1}^{{M_p}} {{{\left| { - \nabla \left( \alpha (x_i) \nabla \widetilde u\left(x_i,  z(x_i)\right) \right) + \beta(x_i)\cdot\nabla \widetilde u\left(x_i,  z(x_i)\right)- f\left( {{{x}_i}} \right)} \right|}^2}}, \\
	{l_2}\left( {\widetilde u} \right) &= \frac{{\left| {{\Omega_s }} \right|}}{{{M_s}}}\sum\limits_{i = 1}^{{M_s}} {{{\left| { - \nabla \left( \alpha (x_i) \nabla \widetilde u\left(x_i,  z(x_i)\right) \right) + \beta(x_i)\cdot\nabla \widetilde u\left(x_i,  z(x_i)\right) - f\left( {{{x}_i}} \right)} \right|}^2}}, \\
	{l_3}\left( {\widetilde u} \right)& =\omega_p \frac{{\left| {{\Gamma }} \right|}}{{{M_{\Gamma}}}}\sum\limits_{i = 1}^{{M_{\Gamma}}} {{{\left| \widetilde u\left(\s^+_i,  z(\s_i^+)\right) - \widetilde u\left(\s^-_i,  z(\s_i^-)\right) - p(\s_i) \right|}^2}}, \\
	{l_4}\left( {\widetilde u} \right) &= \frac{{\left| {{\Gamma }} \right|}}{{{M_{\Gamma}}}}\sum\limits_{i = 1}^{{M_{\Gamma}}} {{{\left|  \alpha(\s_i^+)\frac{\partial u\left(\s^+_i,  z(\s_i^+)\right)}{\partial\nn} - \alpha(\s_i^-)\frac{\partial u\left(\s^-_i,  z(x_i^-)\right) }{\partial\nn} - q(\s_i) \right|}^2}}, \\
	{l_5}\left( {\widetilde u} \right) &= \frac{{\left| {\partial \Omega } \right|}}{{{M_{\partial \Omega }}}}\sum\limits_{i = 1}^{{M_{\partial \Omega }}} { {{{\left| { {\tilde u\left( {{{\s}_i},\phi \left( {{{\s}_i}} \right)} \right) - g\left( {{{\s}_i}} \right)} } \right|}^2}}. } \nonumber
\end{align*}

\subsection{Levenberg-Marquardt optimizer}
To train the parameters in neural networks, some optimizaiton methods must be used to solve the minimization problem \eqref{dlf_noi} or \eqref{dlf_i}. Currently, the Adam optimizer \cite{kingma2014adam} is a popular choice for deep learning methods. However, it seems it limits the accuracy of the method since as shown in \cite{tseng2023cusp,li2023continuity}, using the Levenberg-Marquardt (LM) method as the optimizer, the accuracy using the deep learning method can be significantly improved, the relative errors reaching $10^{-6}$ or even better and potentially outperforming the classic finite element method in terms of solution accuracy. For better performance purpose, we adopt the LM method as the optimizer in this paper. For completeness, we briefly introduce the LM method and its details in our implementation.

For minimization problems encountered in deep learning methods, consider the following general loss function
\[\mbox{Loss}(P) = \frac{1}{M}\sum_{i=1}^M(\phi(x_i) -\phi_{aug}(x_i;P) )^2,\]
where $P\in\mathbb{R}^{N_p}$ denotes the vector of hyperparameters in the specified neural network, and $x_i$ $(i=1,\cdots,M)$ are the sampling points. According to the description in \cite{marquardt1963algorithm}, using the LM method to get the minimal of the above loss function, we have the following iteration scheme:
\[P^{(k+1)} = P^{(k)}+\delta,\]
where $P^{(k)}$ denotes current approximation of $P$, $\delta$ satisfies the following modified normal equation ($\frac{1}{M}$ does not affect the minimization problem)
\[\left( J^TJ+\mu I \right)\delta = J^T(\Phi - \Phi_{aug}(P^{(k)}) ),\]
where $\mu$ is a damping parameter, $\Phi=\{\phi(x_1),\cdots,\phi(x_M) \}$ is the column vector,  
$$\Phi_{aug}(P^{(k)}) =\{ \phi_{aug}(x_1; P^{(k)}),\phi_{aug}(x_2; P^{(k)}), \cdots, \phi_{aug}(x_M; P^{(k)}) \}$$
is another column vector, and $J=\left( \frac{\partial \phi_{aug}(x_i; P^{(k)}) }{\partial P_j}  \right) \in \mathbb{R}^{M\times N_p}$ is the Jacobian matrix.

To use $\frac{1}{M}$ as a scale to possibly improve the condition number of the linear system, we can also rewrite the loss function as follows:
\begin{equation}
\label{gmp}
\mbox{Loss}(P) =\sum_{i=1}^M(\frac{1}{\sqrt{M}}\phi(x_i) -\frac{1}{\sqrt{M}}\phi_{aug}(x_i; P) )^2,
\end{equation}
from which we can similarly derive the scheme to update approximations about $P$. This new scaled scheme will be used in our implementation.

For clarity, we present the LM algorithm in the following:

\begin{minipage}{0.9\linewidth}
\begin{algorithm}[H]
\SetAlgoLined
\hspace{0.5cm}To find the minimal point of the loss function \eqref{gmp}: \\
\textbf{Step 1} : Set initially the tolerance $\epsilon$ and one large value of the damping parameter $\mu$ (by default set $\mu=\mu_{max}=10^8$) so that the steep decent method is used to do the minimization; Fix $\mu_{div} = 3$ and $
\mu_{mul} = 2$; Initialize the hyperparameter values of $P$, loss and others; \\
\textbf{Step 2} : For the first step: use the LM method to update the hyperparameter values of $P$, compute the new loss, the Jacobian matrix $J$, and set $\mu=\max\{\mu/\mu_{div}, 10^{-15}\}$;\\
\textbf{Step 3}:  For other steps: for current approximation $P^{k}$, compute the increment $\delta$ and the new loss with updating hyperparameter values of $P^{k+1}$ using the LM method, then
\begin{minipage}{0.93\linewidth}
\begin{enumerate}
\item If new loss is smaller than the previous loss value, set $\mu=\max\{\mu/\mu_{div}, 10^{-15}\}$ and back to step 3 for the new approximate $P^{k+1}$;
\item Otherwise, compute $\cos\theta$ = $\frac{<\delta,\delta_{old}>}{\| \delta_{old}\| \cdot \| \delta \|}$, the cosine of the angle between previous and current search directions, where $\delta_{old}$ is the increment in previous step: 
\begin{enumerate}
\item If (1-$\cos\theta$) multiplies current loss is great than previous loss (or the minimal of loss have computed so far), then give up current search direction and set $\mu=\max\{\mu*\mu_{mul}, \mu_{max}\}$. With the new $\mu$ value, solve the modified normal equation to get new search direction $\delta$ and back to 3;
\item Otherwise, accept current search direction,  then update necessary information, set $\mu=\max\{\mu/\mu_{div}, 10^{-15}\}$ and back to 3;
\item Stop the iterations when the loss is less than a given tolerance $\epsilon$;
\end{enumerate}
\end{enumerate}
 \end{minipage}
 \caption{ The Levenberg-Marquardt Optimizer}
\end{algorithm}
 \end{minipage}

\subsection{Adaptive failure informed training points design}
\label{set:FI-PINN}
The inclusion of training points is imperative in the development of deep learning-based PDE solvers, especially crucial for training PINNs. With respect to this point,  the FI-PINNs proposed in \cite{gao2023failure} introduces an innovative failure probability estimation framework addressing the challenge of designing optimal training points for PINNs. Inspired by the FI-PINNs proposed and other pertinent adaptive sampling methodologies \cite{lu2021deepxde},  we try to develop a better method to adaptively allocate the sampling points.

First of all, we provide a concise framework for estimating the failure probability in the context of FI-PINNs and here adhere to the mathematical symbols presented in \cite{gao2023failure} for simplifying notations purpose.

The limit-state function (LSF) $g(x):\Omega\to \mathbb{R}$ of the PINNs 
is defined as:
\begin{equation}\label{Power_function}
g(x) = |r(x;\mb{\theta})| - \epsilon_{r},
\end{equation}
where $r(x;\mb{\theta})$ is the residual, $\mb{\theta}$ denotes all hyperparameters of the neural network, and $\epsilon_r$ is a pre-defined parameter.
The failure hypersurface defined by $g(x)=0$ divides the physical domain into two subsets: the safe set $\Omega_{\mathcal{S}} = \{x:g(x) < 0\}$ and the failure set $\Omega_{\mathcal{F}} = \{x:g(x) > 0\}.$   
To describe the reliability of the PINNs, the failure probability $P_{\mathcal{F}}$ is defined under a prior distribution $\omega(x)$:
 \begin{equation}\label{failure_probability}
P_{\mathcal{F}} = \int_{\Omega}\omega(x)\mathbb{I}_{\Omega_{\mathcal{F}}}(x)dx,
\end{equation}
where $\mathbb{I}_{\Omega_{\mathcal{F}}} :\Omega \rightarrow \{0,1\}$ represents the indicator function, taking $\mathbb{I}_{\Omega_{\mathcal{F}}}(x)=1$ when $x \in \Omega_{\mathcal{F}}$, and $\mathbb{I}_{\Omega_{\mathcal{F}}}(x)=0$ otherwise.

The PINNs are deemed reliable when the failure probability is below a specified tolerance, denoted as $\epsilon_{p}$. Conversely, if the failure probability exceeds this threshold, it implies an unreliable system, necessitating enhancements in performance of the PINNs. In essence, a diminishing failure region corresponds to a reduced failure probability, indicating an improved reliable system. This insight motivates the development of adaptive strategies to gradually enhance the reliability of the PINNs, wherein new training points are strategically added into the failure region $\Omega_{\mathcal{F}}$ and the PINNs are subsequently retrained until $P_{\mathcal{F}}\leq \epsilon_{p}$.

To estimate the failure probability, a natural approach involves leveraging MC methods. In this context, one can generate a set of randomly distributed locations $\mathcal{S} = \left\{\mathbf{x}_{1}, \mathbf{x}_{2},\ldots, \mathbf{x}_{|\mathcal{S}|}\right\}$ based on the prior distribution $\omega(\mathbf{x})$ (e.g., the uniform distribution). The MC estimator is then expressed as follows:
$$\hat{P}_{\mathcal{F}}^{MC} = \frac{1}{|\mathcal{S}|}\sum_{x\in \mathcal{S}}\mathbb{I}_{\Omega_{\mathcal{F}}}(x).$$ 
Using this estimator, it becomes possible to formulate an adaptive sampling scheme. If $\hat{P}_{\mathcal{F}}^{MC}>\epsilon_{r}$, new collocation points meeting the criteria $\{\mathbf{x}_{i}:\mathbf{x}_i \in \mathcal{S}, g(\mathbf{x}_{i}) > 0\}$ are incorporated into the training set, followed by a network retraining.
It is worth mentioning that the RAR method \cite{lu2021deepxde} simply selects $m$ new points exhibiting the largest residuals within $\mathcal{S}$.

Given that the failure region might be relatively small compared to the entire domain, the previously mentioned MC sampling strategy may be ineffective in producing useful samples. 
To address this limitation, one can turn to the importance sampling (IS), which has the potential to generate effective samples by mitigating variance through the selection of suitable proposal distributions. In this scenario, the estimation of the failure probability can be achieved as follows:
\begin{equation}
         \label{failure_probabilty_importance_sampling}
         P_{\mathcal{F}} = \int_{\Omega}\mathbb{I}_{\Omega_{\mathcal{F}}}(x)\frac{\omega(x)}{h(x)}h(x)dx =
          \mathbb{E}_{h}\left[\mathbb{I}_{\Omega_{\mathcal{F}}}(x)R(x)\right],
\end{equation}
where $h(x)$ is a proposal distribution, and $R(x)$ represents the weight function $\frac{\omega(x)}{h(x)}$ that transfers the proposal distribution $h(x)$ to the prior distribution $\omega(x)$.  By choosing a set of randomly distributed samples $\mathcal{S} = \{x_{1}, x_{2}, \ldots, x_{|\mathcal{S}|}\}$ from the proposal distribution $h$, we then can similarly approximate the failure probability $P_{\mathcal{F}}$ as
     \begin{equation}
         \label{3.1.3}
\hat{P}_{\mathcal{F}}^{IS} = \frac{1}{|\mathcal{S}|}\sum_{x\in \mathcal{S}}\mathbb{I}_{\Omega_{\mathcal{F}}}(x)R(x).
     \end{equation}
If the support of the proposal distribution $h$ contains the intersection of the support of $\omega$ and the failure set, the equation \eqref{3.1.3} gives an unbiased estimator for $P_{\mathcal{F}}$.  Theoretically, there exists an optimal proposal distribution,
 \begin{equation}
    \label{optimal_importance_distribution}
    h_{opt}(x) = \frac{\mathbb{I}_{\Omega_{\mathcal{F}}} (x)\omega(x)}{P_{\mathcal{F}}} \propto \mathbb{I}_{g(x) > 0}\omega(x),
 \end{equation}
which leads to a zero-variance estimator.  
However, the practical utilization of $h_{opt}(x)$ is impeded by the absence of the normalizing constant. Approximation and adaptive schemes emerge as two widely-used techniques, as evidenced in the literature (\cite{gao2023failure,lu2021deepxde,tang2021deep}), aimed at achieving accurate solutions. Through the adaptive scheme, an approximated optimal proposal distribution $h^i_{opt}(x)$ can be derived for each iteration, with $i$ ranging from 1 to $N_{d}$. Subsequently, the training collocation points are sampled from the mixture distribution defined as $\omega_{train}(x) = \lambda^0*h^0_{opt}(x) + \dots + \lambda^{N_d}*h^{N_d}{opt}(x)$, where $\lambda^i$ represents the weights, and $h^i_{opt}(x)$ signifies the approximation of the optimal proposal density.

In the failure probability estimation framework introduced by the FI-PINNs \cite{gao2023failure}, a truncated Gaussian-based approximation is employed to generate new training points. This approach leverages the efficiency and ease of calculation associated with the density estimation of a Gaussian distribution.
In Section \ref{set:resample_tech}, we present a novel technique for efficiently approximating the optimal proposal distribution. This involves the use of a piece-wise uniform approximation facilitated by the particle transition which speeds up the sampling. 

\subsection{Piece-wise uniform approximation (PUA) of the optimal proposal distribution}
\label{set:resample_tech}

We observe that $\omega(x)$ actually follows a uniform distribution with a constant PDF value within the loss functions \eqref{cminnointerface} and \eqref{cmin}. The indicator $\mathbb{I}_{g(\mathbf{x})>0}$ delineates the input domain featuring a non-zero PDF value within equation \eqref{optimal_importance_distribution}. While the exact expression for $h_{opt}$ remains intractable, it is discernible that $h_{opt}$ conforms to a step-wise uniform distribution. With these observations, we then can introduce a novel piece-wise uniform density distribution as an approximation for $h_{opt}$ to enhance the training efficiency of the PINNs.
Different from \citep{gao2023failure}, where the specific density expression of $h_{opt}$ is employed to efficiently compute $\hat{P}^{IS}_{\mathcal{F}}$, we adopt a MC method for calculating the failure probability due to the cheap calls of neural network.
This approach has the advantage of circumventing the need of the explicit density expression of $h_{opt}$. Consequently, our focus shifts to develop an effective sampling method for the optimal proposal distribution.

A straightforward approach for obtaining samples of \(h_{\text{opt}}\) is through accepted-rejection sampling methods, wherein a candidate is drawn from the prior distribution and accepted if it falls within the failure domain. However, this method becomes inefficient when the failure domain is relatively small compared to the whole domain \(\Omega\). To overcome this issue, motivated by the Adaptive Bayesian Computation-Sequential Monte Carlo (ABC-SMC) method \cite{del2012adaptive}, one can construct a kernel function \(K({\cdot})\) to transform samples from the density function \(\mu(x)\) to another density \(\nu({x})\). By combining accepted-rejection sampling with a transition kernel, we can easily generate enough number of samples from the optimal density.

\begin{figure}
	\center
	\includegraphics*[width=12cm]{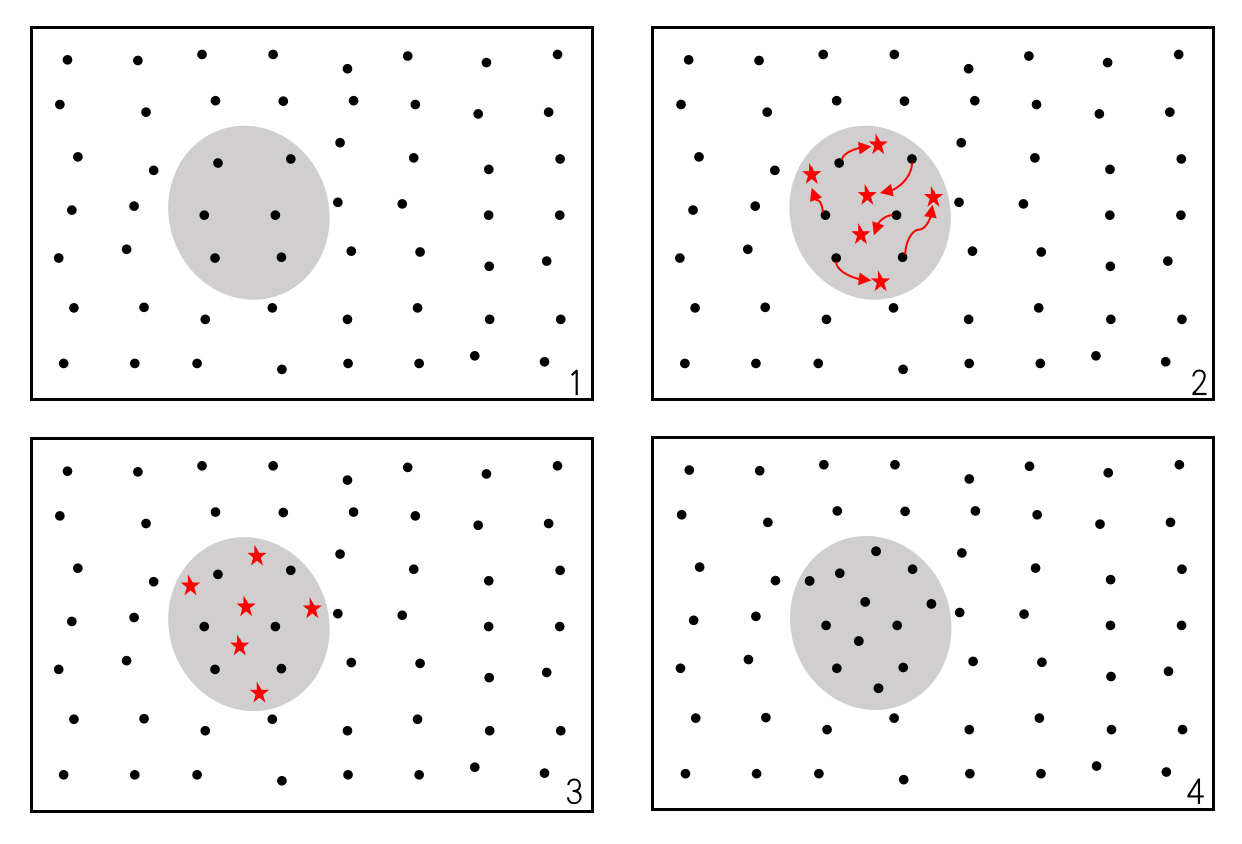}
	\centering
	\caption{Schematic diagram of PUA. The gray area denotes the failure domain, and the black points within this domain indicate the locations where PINNs fail. The red stars represent new sampling points obtained by PUA from the optimal proposal distribution generated by the transition kernel $K(\cdot)$.}
	\label{fig:illustration_PUA}
\end{figure}

The core idea is to employ the accepted-rejection sampling method to obtain initial samples that encompass the entire disconnected failure domain as comprehensively as possible. Subsequently, additional samples are drawn based on the initial set, using a random transition. The transition kernel, denoted as \(K_{\epsilon}(x)\) , has a center at \({x}\) and a random transition length-scale of \(\mathbf{\epsilon}\).
The schematic diagram of this proposed approach, referred to as PUA, is illustrated in Figure \ref{fig:illustration_PUA}.

The transition kernel plays a crucial role in ABC-SMC methodology due to the distinct density expressions of $\mu$ and $\nu$. In this scenario, when $\mu$ and $\nu = \mu$ represent the same distribution, our objective is to efficiently draw more samples from the optimal proposal distribution. To achieve this, a simple and broadly applicable design for the kernel expression, denoted as $K_{\mb{\epsilon}}(\cdot)$, is introduced in the following.

We define the kernel function $K_{\mb{\epsilon}}(x)$ using combinations of the indicator function and  the uniform distribution with means $x$ and step-size $\mb{\epsilon}$ as follows: 
\begin{equation}
\label{eq:transition_kernel}
K_{\epsilon}(x) = x'\otimes \mathbb{I}_{\Omega_{\mathcal{F}}}(x') + x\otimes \mathbb{I}_{\Omega_{\mathcal{S}}}(x'), \quad x'\sim \mathcal{U}(x_i-\mb{\epsilon}, x_i+\mb{\epsilon}),
\end{equation}
where $\otimes$ represents the Kronecker product, ${x}'$ is suggested through a proposal distribution $\mathcal{U}({x}_i - \mb{\epsilon}, {x}_i + \mb{\epsilon})$ to align with the target distribution, given that the current distribution is identical to the target.
The acceptance of the proposed sample hinges on its location within the failure domain, while rejection occurs if the proposed point falls outside the bounds of the failure domain.
Here, the step-size $\mb{\epsilon}$ is computed as 
\begin{equation}
\label{eq:epsilon_def}
\mb{\epsilon} = \frac{|\Omega|}{2|\mathcal{S}|},
\end{equation}
where $|\mathcal{S}|$ is the total number of points used in the accepted-rejection sampling method for $M$ failure points. The intuition behind this calculation is quite obvious, as $2*\mb{\epsilon}$ represents the average area controlled by each training point. For clarity, the detailed PUA algorithm is also listed in the following\\
\begin{algorithm}[H]
    \caption{Piece-wise approximation of the optimal proposal distribution}
    \label{alg:optimal_proposal_sample}
    \begin{algorithmic}[1]
    \STATE{\textbf{Inputs}:  $N_i$ sampling points in the $i$-th iteration, the limit state function $g(x)$, the uniform distribution $\omega(x)$;}
    \STATE{\textbf{Outputs}: Samples set $\mathcal{D}'$ of the optimal proposal;}
    \STATE{Find $M$ failure points $\mathcal{D}^* = \{x^*_{1},\dots,x^*_{M}\}$ by $g(x)$ and $\omega(x)$ by the accepted-rejection sampling method;}
    \STATE Compute the $\mb{\epsilon}$ by the equation \eqref{eq:epsilon_def} and set $j=1$;
\WHILE {$j \le N_i$}	
	\STATE Randomly draw $x^*$ from $\mathcal{D}^*$;
	\STATE Generate new points $x_j'$ based on $x^*$ and $K_{\epsilon}(\cdot)$ given in the equation \eqref{eq:transition_kernel};
	\STATE $j=j+1$;
\ENDWHILE
\STATE {$\mathcal{D}' = \{ x_1',\dots,x_{N_i}'\}$};
   \end{algorithmic}
\end{algorithm}

\section{Numerical results}
In this section, we carry out the numerical experiments using the methods given before. In all examples, we choose the activation function to be hyperbolic function $\tanh(x)$. During each iteration, the program will be terminated when the loss function reaches $1\times 10 ^{-12}$ or the epoch number reaches the maximal number 10000 by default. At the beginning, we always use $M=2000$ collocation points and $M_b=200$ boundary points except the elliptic interface problem to train the network. Without explicit statement, we will employ the PUA method to obtain new points and resample, ensuring that the number of points before and after each iteration remains the same. The experiments are conducted using PyTorch on a 2.2GHz processor (Intel(R) Core(TM) i7-8750H CPU) with 8GB of RAM and Windows OS (Win11).

Meanwhile, in order to testify the method, we use the following relative $L^2$ error
\begin{equation}
      err_{L^{2}}=\frac{\sqrt{\sum_{i=1}^{N}\left | \hat{u}({x}_{i})-u({x}_{i}) \right |^{2} }}{\sqrt{\sum_{i=1}^{N}\left | u({x}_{i}) \right |^{2} }},
\end{equation}
where $N$ represents the total number of test points chosen uniformly in the domain, and $\hat{u}({x}_{i})$ and $u({x}_{i})$ represent the predicted and the exact solution values at point $x_i$, respectively. For the following test cases, we set $N=500^2$, $N=30^3$, and $N=20^5$ for the two dimensional, the three dimensional, and the five dimensional test cases, respectively.

\subsection{Two-dimensional Poisson equation}\label{example1}
We first consider the following two-dimensional Poisson equation
\begin{align*}
\left\{\begin{array}{ll}
            -\Delta u(x,y)=f(x,y)\quad  &\mbox{in } \Omega,\\ 
             u(x,y)=g(x,y) \quad &\mbox{on }\partial \Omega,
\end{array}\right.
\end{align*}
where $\Omega$ is $[-1,1]^{2}$ and the exact solution is taken to be as follows
\begin{equation}
      u(x,y)=\exp(-1000[(x-0.5)^{2}+(y-0.5)^{2}])\;,
\end{equation}
which has a peak at (0.5, 0.5) and decreases rapidly away from (0.5, 0.5).

\begin{figure}[h]
    \centering
    \includegraphics[width=0.35\linewidth]{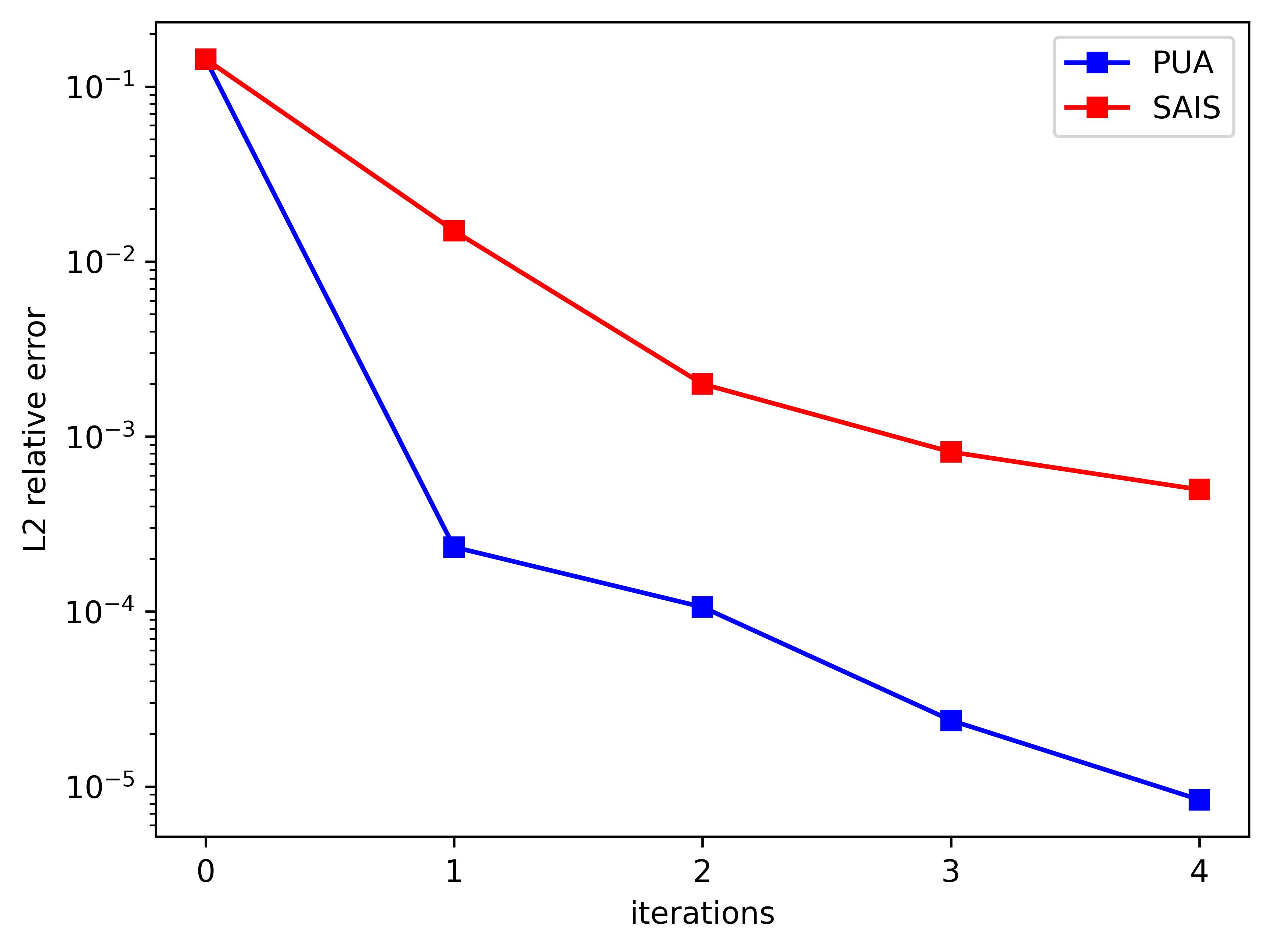}
    \caption{The relative $L^2$ errors for example \ref{example1} obtained by two different adaptive sampling methods.}
    \label{fig:example1_l2}
\end{figure} 

\begin{figure}[h]
    \centering
    \includegraphics[width=6in]{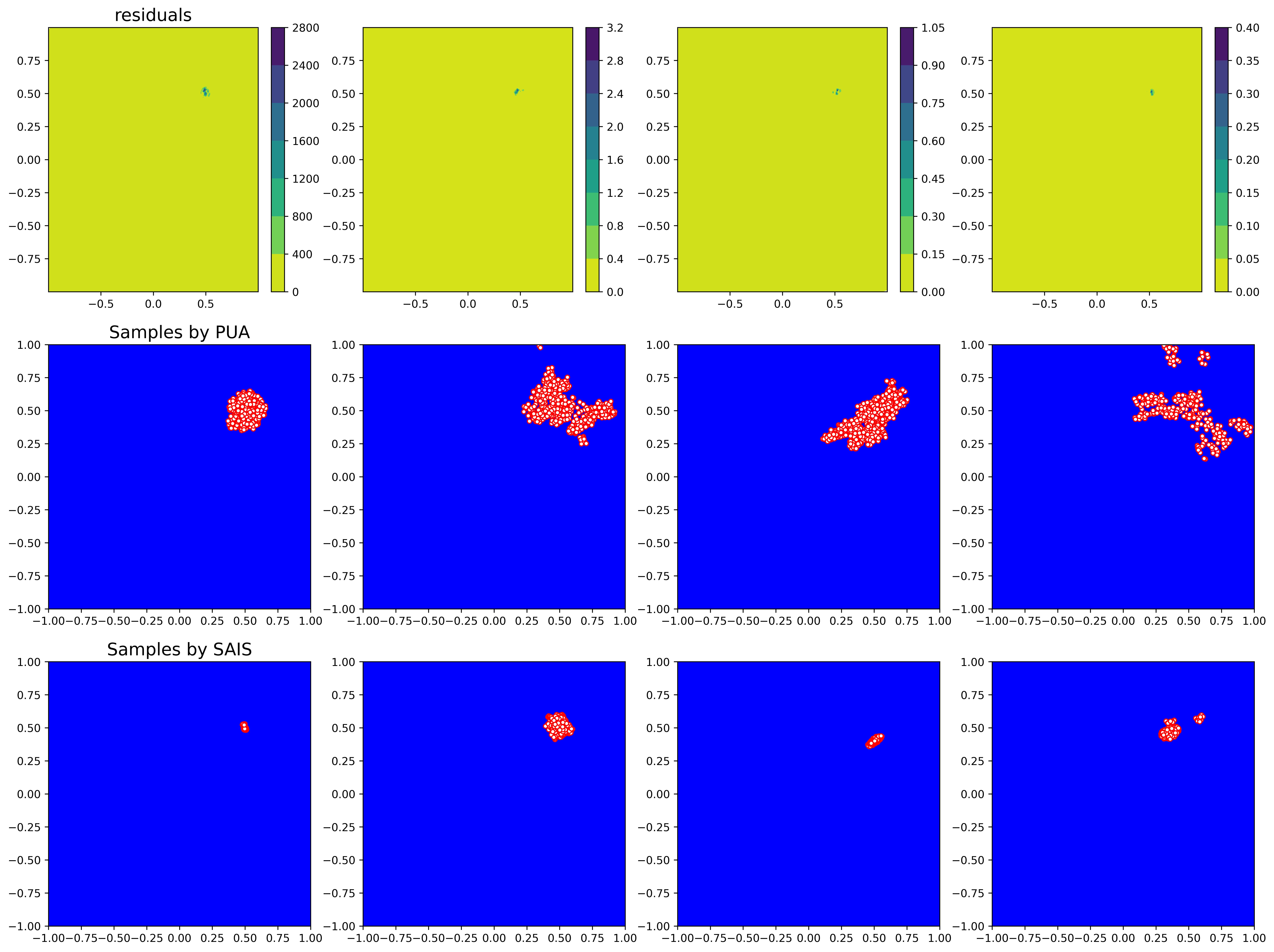}
    \caption{Contour plots of residuals of the loss function and distributions of new points obtained by PUA and SAIS for four iterations.}
    \label{fig:example1_samples}
\end{figure}

\begin{figure}[!ht]
    \centering
    \includegraphics[width=6in]{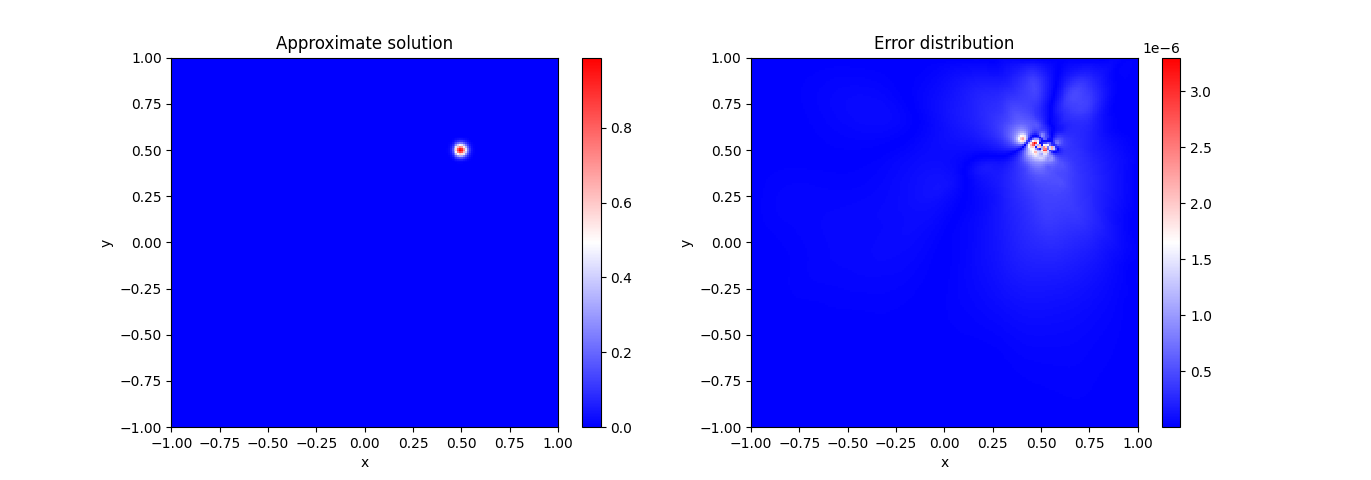}
    \caption{Left: The solution profile obtained by PUA. Right: Profile of the pointwise errors. Here the relative $L^2$ error is $8.41\times 10^{-6}$.}
    \label{fig:example1_pred}
\end{figure}

For this example, we testify the improvements of our new method compared to the SAIS method proposed in \cite{gao2023failure}. To carry out the computation, the structure of the neural network is taken to be $(L, N_L, N_p) = (7, 20, 2600)$, where $L$ is the number of the hidden layers, $N_L$ is the number of neurons per hidden layer and $N_p$ is the total number of hyperparameters in the neural network. 
To make a fair comparison, we firstly fully train the network using an uniformly distributed dataset, and then use two different adaptive methods (The PUA and the SAIS methods) to compare the performance. 
Here not only we ensure that the threshold $\epsilon$ for the PUA method and the residual tolerance $\epsilon_{r}$ used in SAIS are equal, but also the generated numbers of new points using the PUA method and the SAIS method are exactly the same. The parameters we use, the corresponding numbers of new points, and the obtained relative $L^2$ errors are shown in Table \ref{table:example1}. The comparison results are presented in Fig \ref{fig:example1_l2} and Fig \ref{fig:example1_samples} , while the obtained solution as well as the error distributions after 4 iterations of the PUA method are given in Fig \ref{fig:example1_pred}. It can be seen that the error achieved by the PUA method decreases much faster than the SAIS method,
and the relative $L^2$ error is much smaller, which is $8.41\times 10 ^{-6}$ compared to $4.98\times 10 ^{-4}$ obtained by SAIS. From these numerical results, it confirms that the PUA method is indeed an improved variant of the SAIS method. Additionally, we observed that the time required for performing PUA can be nearly equivalent to that of SAIS when employing a single Gaussian distribution to approximate $h_{opt}$.


\begin{table}[ht]
    \centering
    \begin{tabular}{l|ccccc}
        \hline
        Iterations & 0 & 1 & 2 & 3 & 4 \\
        \hline
        $\epsilon$ ($\epsilon_{r} $) & & 0.1 & 0.001 & 0.0005 & 0.0001 \\
        Number of new points & & 641 & 812 & 993 & 375 \\
         Rel $L^2$ errors (PUA) & $1.44\times 10 ^{-1}$ & $2.34\times 10 ^{-4}$ & $1.06\times 10 ^{-4}$ & $2.39\times 10 ^{-5}$ & $8.41\times 10 ^{-6}$ \\
        Rel  $L^2$ errors (SAIS)& $1.44\times 10 ^{-1}$ & $1.50\times 10 ^{-2}$ & $2.00\times 10 ^{-3}$ & $8.17\times 10 ^{-4}$ & $4.98\times 10 ^{-4}$\\
        \hline
    \end{tabular}
    \caption{Parameters and numbers of new points as well as the errors for example \ref{example1}.}
    \label{table:example1}
\end{table}

\subsection{Two-dimensional problem with 2 peaks}\label{example2}
Here we consider the following two-dimensional PDE, solution of which contains 2 peaks,
	\begin{align*}
	\left\{\begin{array}{ll}
		-\triangledown \cdot [u(x,y)\bigtriangledown (x^{2}+y^{2})]+\bigtriangledown ^{2}u(x,y)=f(x,y) \quad &\mbox{in}\; \Omega,\\ 
		   u(x,y)=b(x,y) \quad &\mbox{on}\; \partial \Omega,
       \end{array}\right.
	\end{align*}
where $\Omega$ is $[-1,1]^{2}$. The true solution is taken to be 
\begin{equation}
	u(x,y)=\sum_{i=1}^{2}\exp(-1000[(x-x_{i})^{2}+(y-y_{i})^{2}])\;,
\end{equation}
where $(x_1,y_1)=(0.5, 0.5)$ and $(x_2,y_2)=(-0.5, -0.5)$. We can easily see the solution of the problem has peaks at $(0.5, 0.5)$ and $(-0.5, -0.5)$, and decays to zero exponentially in other places.

\begin{table}[!ht]
    \centering
    \begin{tabular}{l|ccccc}
        \hline
        Iterations & 0 & 1 & 2 & 3 & 4 \\
        \hline
        $\epsilon $ & & 1 & 0.01 & 0.005 & 0.001 \\
        Number of new points & & 824 & 1106 & 541 & 995 \\
        Rel $L^2$ errors & $4.40\times 10 ^{-1}$ & $4.60\times 10 ^{-4}$ & $4.20\times 10 ^{-5}$ & $2.98\times 10 ^{-5}$ & $1.98\times 10 ^{-5}$ \\
        \hline
    \end{tabular}
    \caption{Parameters and numbers of new points as well as the errors for example \ref{example2} with setting $(M, M_b)= (2000, 200)$. }
    \label{table:example2_1}
\end{table}
\begin{table}[!ht]
    \centering
    \begin{tabular}{l|ccccc}
        \hline
        Iterations & 0 & 1 & 2 & 3 & 4 \\
        \hline
        $\epsilon $ & & 1 & 0.005 & 0.001 & 0.00025 \\
        Number of new points & & 1206 & 882 & 1118 & 1138 \\
        Rel  $L^2$ errors & $1.08\times 10 ^{-1}$ & $3.90\times 10 ^{-5}$ & $8.50\times 10 ^{-6}$ & $4.86\times 10 ^{-6}$ & $3.25\times 10 ^{-6}$ \\
        \hline
    \end{tabular}
    \caption{Parameters and numbers of new points as well as the errors for example \ref{example2} with setting $(M, M_b) = (5000, 500)$.} 
    \label{table:example2_2}
\end{table}

\begin{figure}[!ht]
    \centering
    \includegraphics[width=6in]{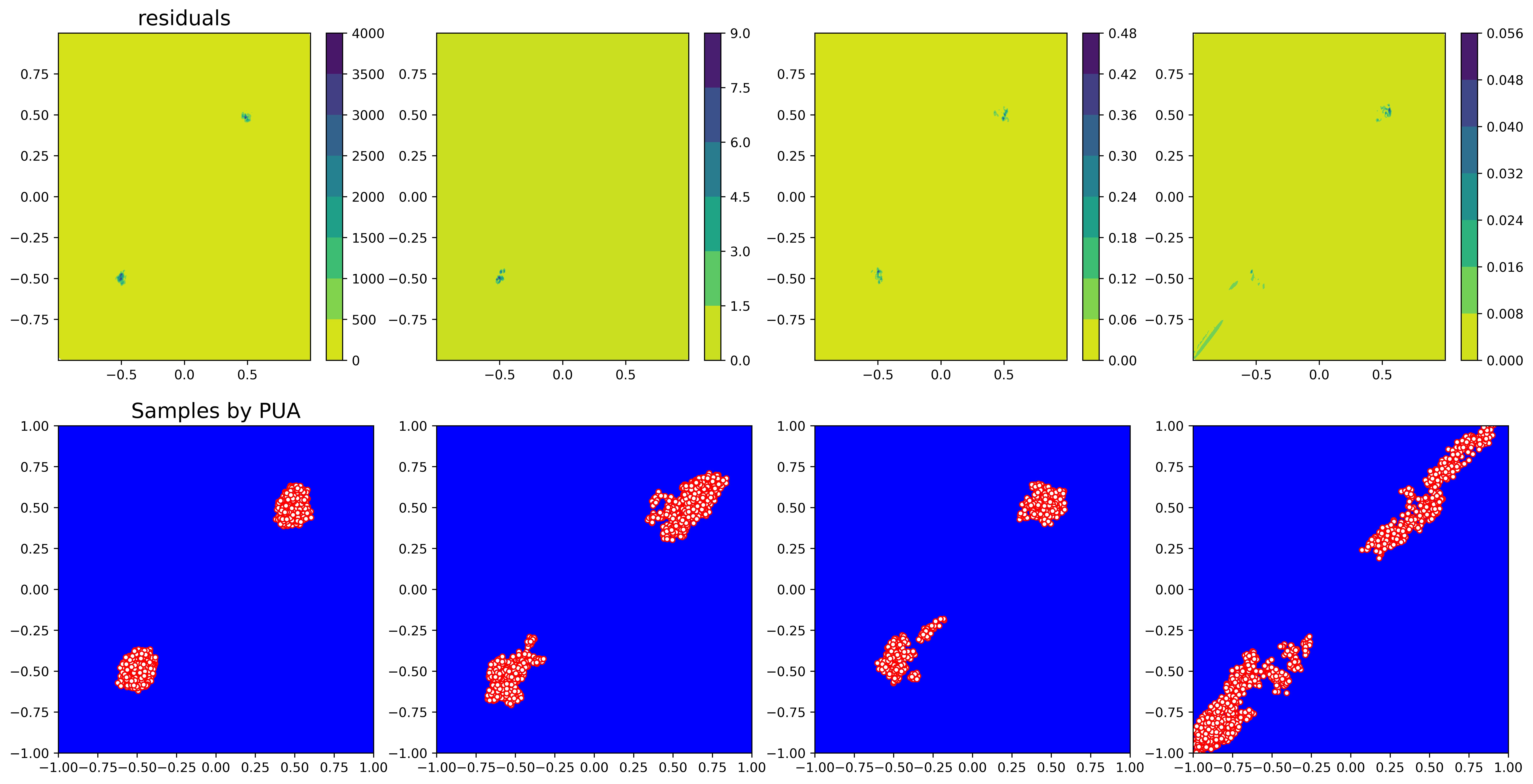}
    \caption{Contour plots of residuals of the loss function and distributions of new points obtained by PUA for four iterations.}
    \label{fig:example2_samples}
\end{figure}

\begin{figure}[!ht]
    \centering
    \includegraphics[width=6in]{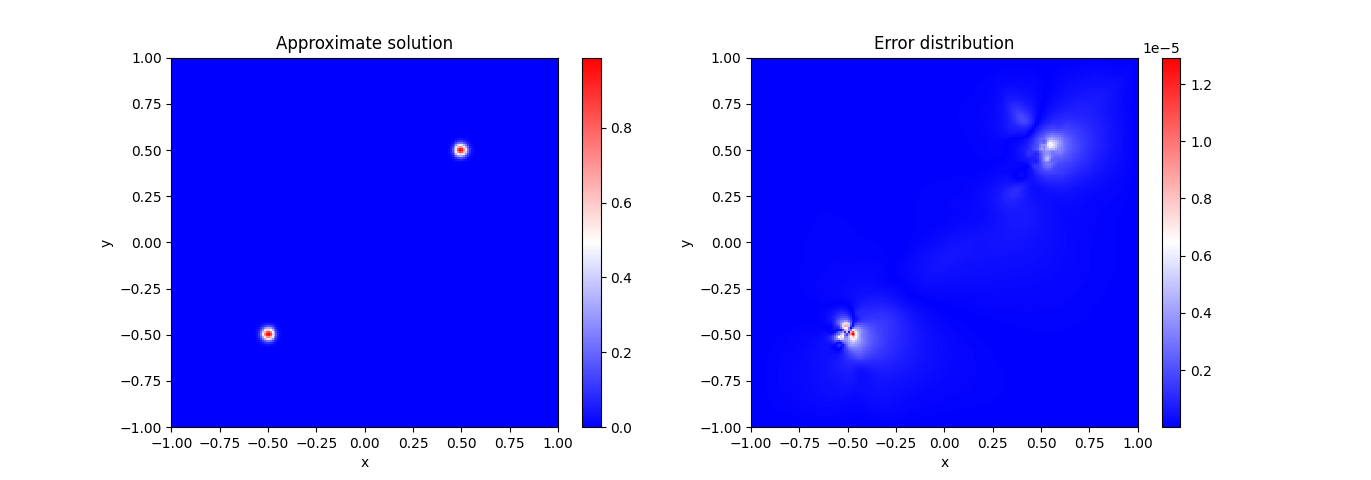}
    \caption{Left: The solution profile obtained by PUA. Right: Profile of the pointwise errors. Here the relative $L^2$ error is $1.98\times 10 ^{-5}$.}
    \label{fig:example2_pred}
\end{figure}

The structure of the neural network is again taken to be $(L, N_L, N_p) = (7, 20, 2600)$. And in the PUA method, the candidate new points are uniformly sampled in $\Omega$ with total $50000$ points.
For this case, we use two different sets of $(M, M_b)$ to testify the effectiveness of the proposed adaptive strategy in this paper. Here one is the default value $(M, M_b)=(2000, 200)$ and the other set is $(M, M_b) = (5000, 500)$.
For these two sets, the parameters we use in the PUA method, the corresponding numbers of new points, and the obtained relative $L^2$ errors of the numerical solutions are given in Tables \ref{table:example2_1} and \ref{table:example2_2}, respectively. From the Tables, we can see that the relative $L^2$ error before the adaptive sampling is $4.58\times 10 ^{-1}$ and after four adaptive sampling processes drops to $1.98\times 10 ^{-5}$ for case $(M, M_b)$ = (2000, 200), and
the relative $L^2$ error before the adaptive sampling is $6.36\times 10 ^{-2}$ and after four adaptive sampling processes drops to $3.25\times 10 ^{-6}$ for case $(M, M_b)$ = (5000, 500), both of which demonstrate the effectiveness of the proposed adaptive sampling method. Of course, the relative $L^2$ errors after one adaptive sampling process have reached $4.6\times 10^{-4}$ and $3.9\times 10^{-5}$ for these two cases, which show that the method can reduce the relative errors quickly, main advantage of our method compared to the SAIS method. Fig \ref{fig:example2_samples} illustrates residuals of the loss function and distributions of new points obtained by the PUA method for four iterations when $(M, M_b)$ = (2000, 200). Meanwhile, the obtained solution as well as the error distributions after four iterations of the PUA method are given in Fig \ref{fig:example2_pred}.


\subsection{High-dimensional Poisson equation}\label{example3}
Consider the following high dimensional elliptic problem:
	\begin{align*}
	\left\{\begin{array}{ll}
		  -\Delta u(x)=f(x) \quad &\mbox{in}\; \Omega,\\ 
		   u(x)=g(x) \quad & \mbox{on}\; \partial \Omega,
        \end{array}\right.
	\end{align*}
where $\Omega$ is $[-1,1]^{d}$ with an exact solution given by the following
\begin{equation}
	u(\mathbf{x}) = \exp(-10\sum_{i=1}^{d}(x_{i}-0.5)^{2})\;.
\end{equation}
We set $d$ = 5 for this test case. To testify the method for this high dimensional problem, the structure of the neural network is still taken to be $(L, N_L, N_p) = (7, 20, 2660)$ and $(M, M_b)$ is set to be $(5000, 500)$ for the initial training. Meanwhile, in the PUA method, the candidate new points are uniformly sampled from $[-1,1]^{5}$ with total $100000$ points. The parameters we use in the PUA method, the corresponding numbers of new points, and the obtained relative errors are given in Table \ref{table:example3}. From the table, we see that after one iteration of the adaptive sampling method, the relative $L^2$ error decreases from  $6.93\times 10 ^{-4}$ to $2.30\times 10 ^{-6}$, and more iterations lead to further decrease in relative errors, reaching $4.08\times 10 ^{-7}$ after another two iterations. The numerical results show the effectiveness of the proposed method for the high dimensional problems.

\begin{table}[!ht]
    \centering
    \begin{tabular}{l|cccc}
        \hline
        Iterations & 0 & 1 & 2 & 3  \\
        \hline
        $\epsilon $ & & 0.0001 & 0.000005 & 0.000003 \\
        Numbers of new points & & 2300 & 716 & 1606 \\
        Rel $L^2$ errors & $6.93\times 10 ^{-4}$ & $2.30\times 10 ^{-6}$ & $5.53\times 10 ^{-7}$ & $4.08\times 10 ^{-7}$  \\
        \hline
    \end{tabular}
    \caption{Parameters and numbers of new points as well as the errors for example \ref{example3}.} 
    \label{table:example3}
\end{table}

\subsection{Elliptic interface problem}\label{example6}
In this example, we test our method for an elliptic interface problem with high contrast diffusion coefficients. We choose a square domain $\Omega=[-1,1]^{2}$ and the interface to be the zero-level set of the function $\phi(x, y) = x^{2}+y^{2}-0.5^{2}$. The exact solution is chosen as:
\begin{equation}
	u(x,y) = 
	\left\{\begin{array}{cl}
		\frac{r^{3}}{\alpha _{1}},& \quad  \mbox{in }  \Omega_p,\\
		\frac{r^{3}}{\alpha _{2}}+\left ( \frac{1}{\alpha _{1}}-\frac{1}{\alpha _{2}} \right )r_{0}^{3},& \quad \mbox{in } \Omega_s,
	\end{array}\right.
\end{equation}
where $\alpha _{1}=1$, $\alpha _{2}=1000$ and $r=x^{2}+y^{2}$ with $r_0=0.5$.

\begin{figure}[!ht]
    \centering
    \includegraphics[width=4.8in]{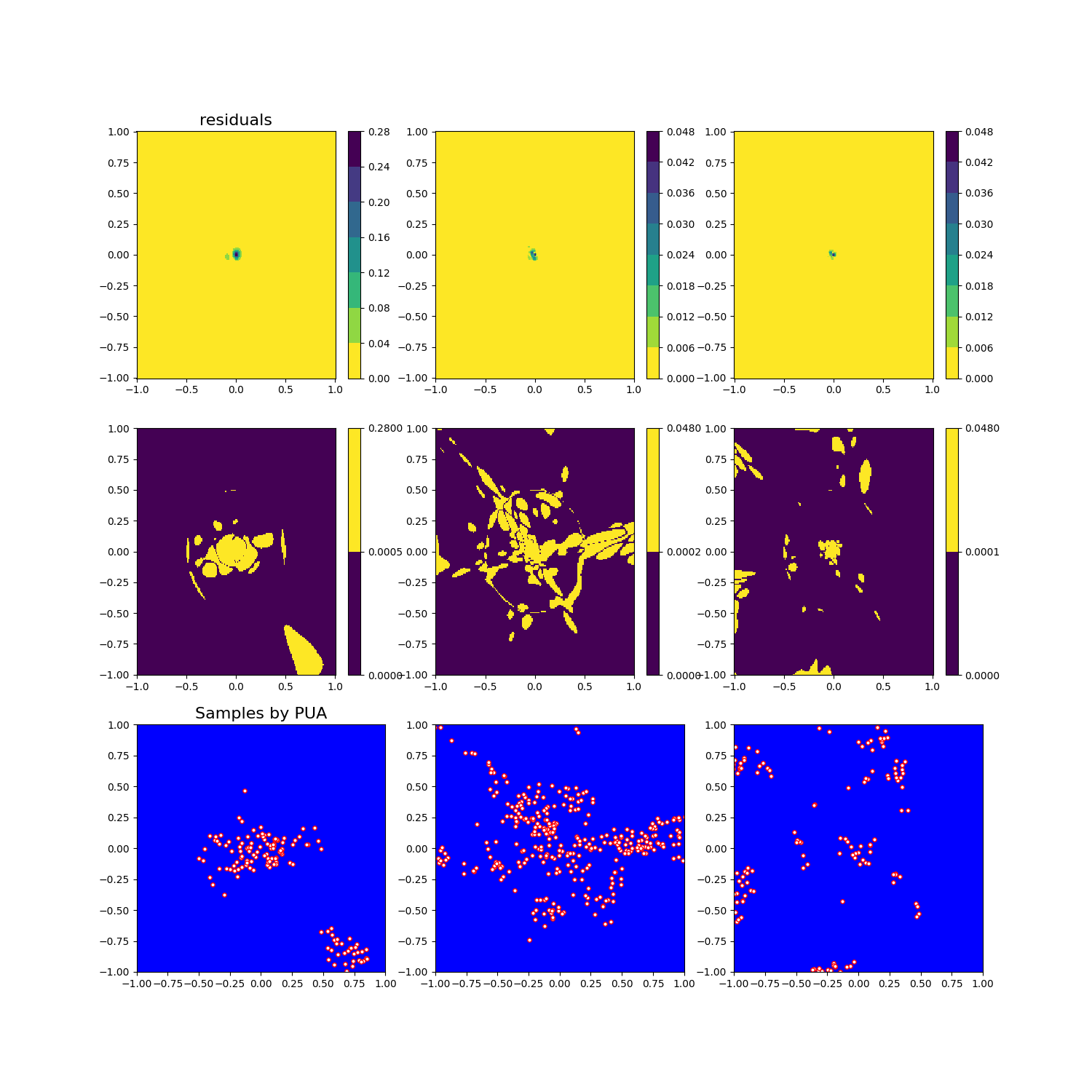}
    \vspace{-1cm}
    \caption{Contour plots of residuals of the loss function and distributions of new points obtained by PUA for three iterations.}
    \label{fig:example6_samples}
\end{figure}

\begin{figure}[!ht]
    \centering
    \includegraphics[width=6in]{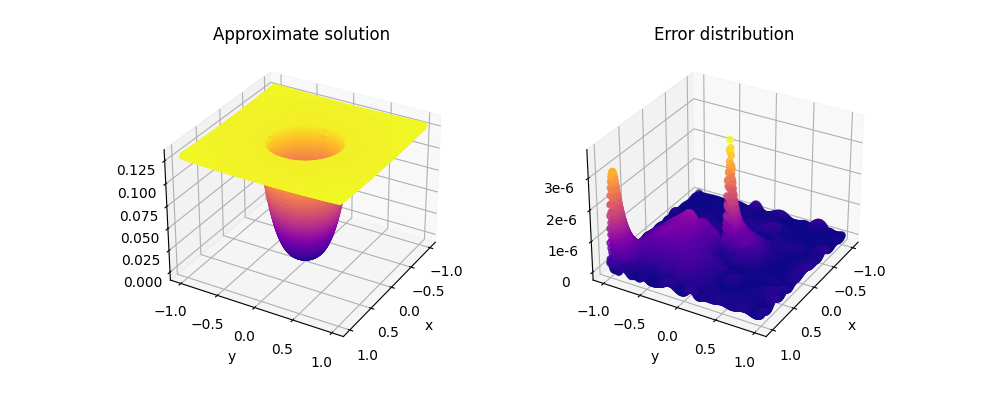}
    \caption{Left: The solution profile obtained by PUA. Right: Profile of the pointwise errors. Here the relative $L^2$ error is $2.99\times 10 ^{-6}$.}
    \label{fig:example6_pred}
\end{figure}

For this interface problem, the structure of the neural network is taken to be $(L, N_L, N_p) = (2, 30, 1080)$ with the initial training points $(M, M_b, M_{\Gamma}) = (500, 100, 100)$, where $M_{\Gamma}$ denotes the number of collocation points on the interface $\Gamma$. 
Different from previous examples, after each iteration in the adaptive sampling process, we add all new points obtained by the PUA method to the current sampling points without resampling to do the retraining. Here in the PUA method, the candidate new points are uniformly sampled from $\Omega$ with total $2500$ points. And the parameters we use in the PUA method, the corresponding numbers of added new points, and the obtained relative errors are given in Table \ref{table:example6}. From the table, we see that the relative $L^2$ error before the adaptive sampling is $4.84\times 10 ^{-4}$ and drops to $2.99\times 10 ^{-6}$ after three adaptive sampling processes. Fig \ref{fig:example6_samples} illustrates residuals of the loss function and distributions of new points obtained by the PUA method for three iterations. Meanwhile, the obtained solution as well as the error distributions after three iterations of the PUA method are given in Fig \ref{fig:example6_pred}, from which we can see that the maximal pointwise error reaches $3\times 10^{-6}$, a quite accurate result already.
\begin{table}[!ht]
    \centering
    \begin{tabular}{l|ccccc}
        \hline
        Iterations & 0 & 1 & 2 & 3 \\
        \hline
        $\epsilon $ & & 0.0005 & 0.0002 & 0.0001 \\
        Numbers of new points($\Omega_1$) & & 34 & 160 & 96 \\
        Numbers of new points($\Omega_2$) & & 99 & 170 & 35 \\
        Rel $L^2$ errors & $4.84\times 10 ^{-4}$ & $5.83\times 10 ^{-5}$ & $4.99\times 10 ^{-6}$ & $2.99\times 10 ^{-6}$  \\
        \hline
    \end{tabular}
    \caption{Parameters and numbers of new points as well as the relative errors for example \ref{example6}.} 
    \label{table:example6}
\end{table}

\subsection{Convection-dominated problem: Case 1}\label{example4}
\begin{figure}[!ht]
    \centering
    \includegraphics[width=6.5in]{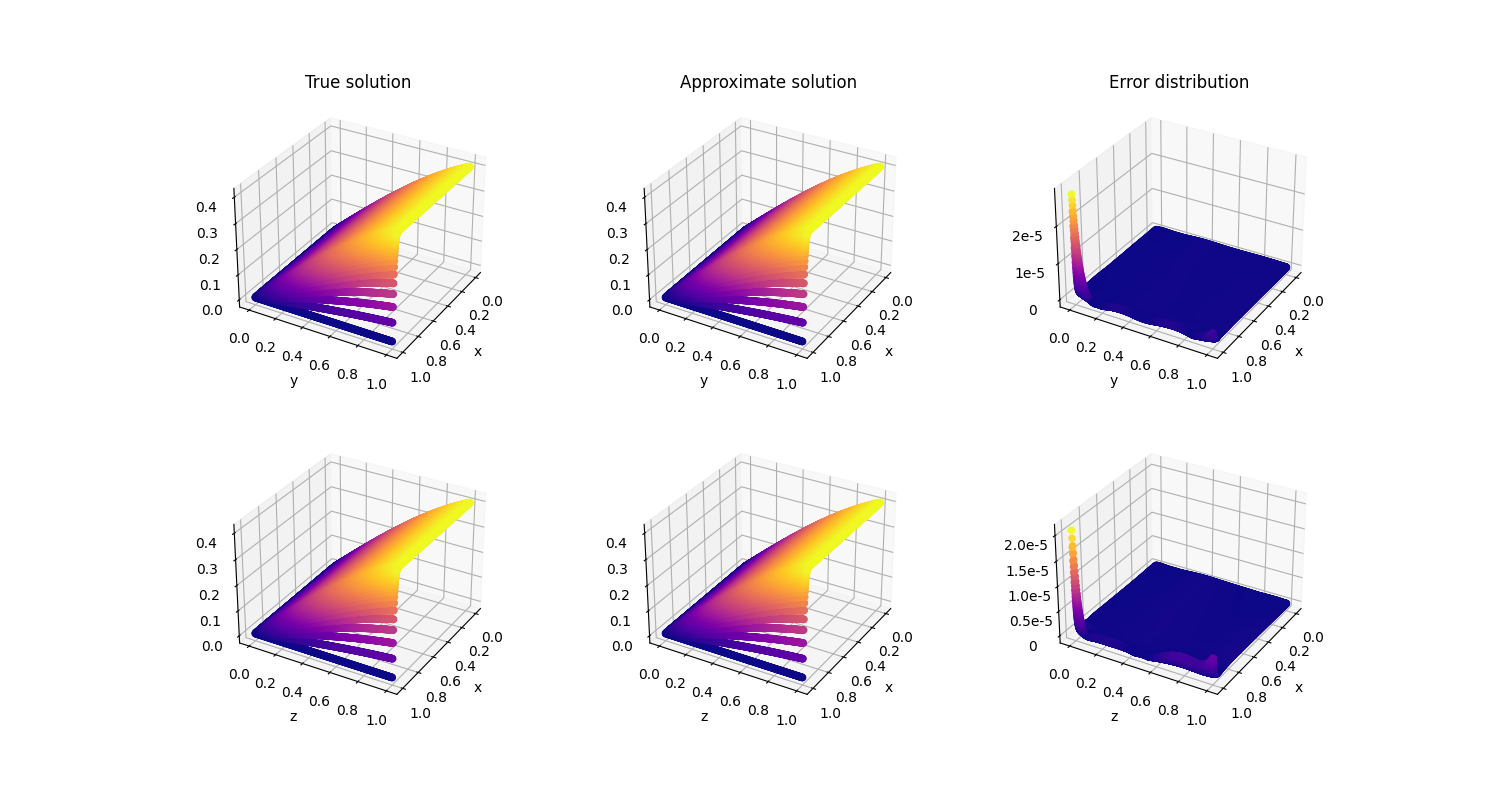}
    \caption{Profiles of the analytical solution and the approximated solution as well as the pointwise errors on slices $z=0.5$ (Top ones) and $y=0.5$ (Bottom ones), respectively.}
    \label{fig:example4_pred}
\end{figure}

This example is from \cite{hus2023efficient}. Consider the 3D steady convection-diffusion equation with constant coefficients:
\begin{equation}
	-(u_{xx}+u_{yy}+u_{zz})+\mathrm{Re}u_{x}=f(x,y,z), \quad \mbox{in }\Omega,
\end{equation}
where $\Omega=[0,1]^3$, and the analytical solution is given by
\begin{equation}
	u(x,y,z)=\frac{(e^{\mathrm{Re}}-e^{\mathrm{Re}x})\sin(y)\sin(z)}{e^{\mathrm{Re}}-1}.
\end{equation}
We set $\mathrm{Re}$ = 100 as done in \cite{hus2023efficient}. And the structure of the neural network is set to be $(L, N_L, N_p)$ = $(3, 30, 2010)$. In the PUA method, the candidate new points are uniformly sampled from $[0,1]^{3}$ with total $50000$ points. The parameters we use, the corresponding numbers of new points, and the numerical results are reported in Table \ref{table:example4}. From the table, we see that after applying the adaptive sampling method three times, the relative $L^2$ error drops from $4.13\times 10 ^{-4}$ to $4.94\times 10 ^{-6}$, significantly improving solution accuracy. It should be pointed out that the relative $L^2$ error obtained in the end actually is better than the one obtained using the $64\times 64\times 64$ nonuniform rectilinear grids reported in \cite[Table 2]{hus2023efficient}. The deep learning method uses far less parameters to obtain better approximation, which is another evidence that the new method can outperform the classic numerical methods in terms of accuracy. Finally, for this 3D test problem, profiles of the analytical solution and the approximated solution as well as the point-wise error distributions on slices $y=0.5$ and $z=0.5$ are given in Fig \ref{fig:example4_pred}, demonstrating that deep learning method is able to capture the boundary layer in the solution.

\begin{table}[!ht]
    \centering
    \begin{tabular}{l|cccc}
        \hline
        Iterations & 0 & 1 & 2 & 3 \\
        \hline
        $\epsilon $ & & 0.5 & 0.01 & 0.002 \\
        Number of new points & & 245 & 215 & 148 \\
        Rel $L^2$ errors & $4.13\times 10 ^{-4}$ & $1.20\times 10 ^{-5}$ & $7.55\times 10 ^{-6}$ & $4.94\times 10 ^{-6}$  \\
        \hline
    \end{tabular}
    \caption{Parameters and numbers of new points as well as the errors for example \ref{example4}.} 
    \label{table:example4}
\end{table}

\subsection{Convection-dominated problem: Case 2}\label{example5}
This example is also from \cite{hus2023efficient}. Consider a 3D convection-diffusion equation with variable coefficients:
\begin{equation}
	-\varepsilon (u_{xx}+u_{yy}+u_{zz})+pu_{x}+qu_{x}+ru_{x}=f(x,y,z),\quad \mbox{in }\Omega, 
\end{equation}
where $\Omega=[0,1]^3$, the functions $p$, $q$, $r$ are defined by
	\begin{align*}
		p=-x(1-y)(2-z),\quad
		q=-(2-x)y(1-z),\quad
		r=-(1-x)(2-y)z,
	\end{align*}
and the analytical solution is taken to be the following function
\begin{equation}
	u(x,y,z)=\frac{e^{\frac{x}{\varepsilon }}+e^{\frac{y}{\varepsilon }}+e^{\frac{z}{\varepsilon }}-2}{e^{\frac{1}{\varepsilon }}-1}.
\end{equation}

\begin{figure}[!ht]
    \centering
    \includegraphics[width=6.5in]{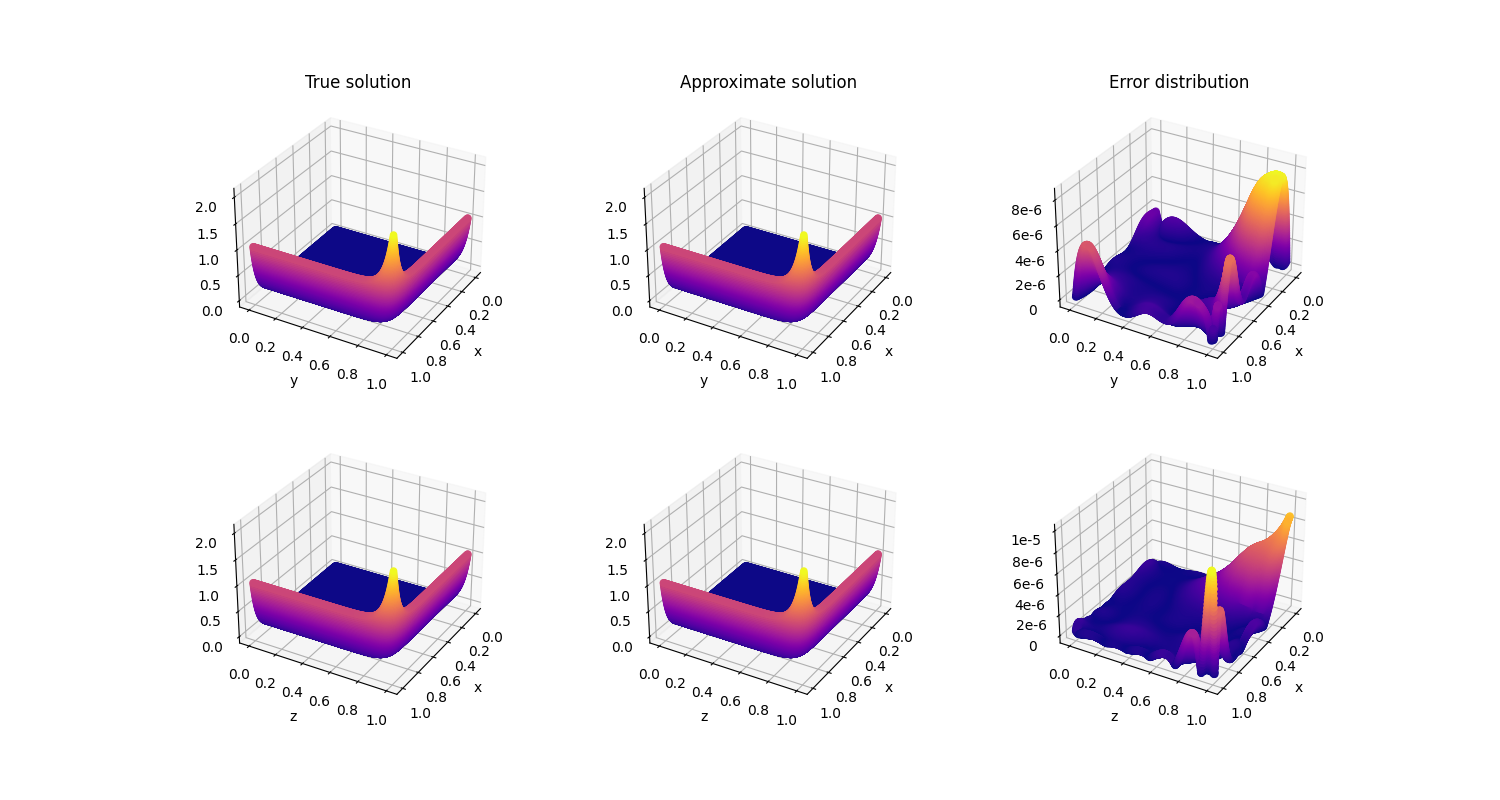}
    \caption{Profiles of the analytical solution and the approximated solution as well as the pointwise errors on slices $z=0.8$ (Top ones) and $y=0.8$ (Bottom ones), respectively.}
    \label{fig:example5_pred}
\end{figure}

We set $\varepsilon = 0.05$ in the equation. For this case, the structure of the neural network is also taken to be $(L, N_L, N_p) = (3, 30, 2010)$, same as the  previous example. In the PUA method, the candidate new points are uniformly sampled from $[0,1]^{3}$ with total $50000$ points. The parameters we use, the corresponding numbers of new points, and the obtained relative errors are reported in Table \ref{table:example5}. From the table, we see that
the relative $L^2$ error before the adaptive sampling is $1.76\times 10 ^{-4}$ and it drops to $3.93\times 10 ^{-6}$ after three adaptive sampling processes.
At last, for the convection-dominated test problem, profiles of the analytical solution and the approximated solution as well as the point-wise error distributions on slices $y=0.8$ and $z=0.8$ are given in Fig \ref{fig:example5_pred}.

\begin{table}[!ht]
    \centering
    \begin{tabular}{l|ccccc}
        \hline
        Iterations & 0 & 1 & 2 & 3 \\
        \hline
        $\epsilon $ & & 0.001 & 0.002 & 0.0001 \\
        new points & & 269 & 251 & 266 \\
        relative L2 & $1.76\times 10 ^{-4}$ & $5.83\times 10 ^{-6}$ & $4.88\times 10 ^{-6}$ & $3.93\times 10 ^{-6}$  \\
        \hline
    \end{tabular}
    \caption{Parameters and numbers of new points as well as the relative errors for example \ref{example5}.} 
    \label{table:example5}
\end{table}


\section{Conclusions}
In this paper, based on the LM method as the solver of the optimization problem, a new and improved adaptive deep learning method is proposed to solve elliptic problems, including the interface problems and the convection-dominated problems. Inspired by the FI-PINNs, a novel approach employs a piece-wise uniform distribution to approximate the optimal proposal distribution for determining training locations optimally. To efficiently draw samples from this distribution, a kernel-based sampling algorithm has been devised. 
The combination of these two components leads to a new and better adaptive method, compared to the SAIS method used in \cite{gao2023failure}.
Numerical tests on various examples confirm this claim. Considering the elliptic problems with or without interface conditions and the convection-dominated problems, numerical results show that by the new adaptive method,  the relative errors not only reduce much faster than the SAIS method, but also are far more better than the ones obtained without any using adaptive sampling technique, the reducing factors of the relative errors reaching $10^4$ for best cases.

\section*{Acknowledgments}
The first author was financially supported by the Natural Science Foundation of Hunan Province (Grant No. 2023JJ30648) and the Natural Science Foundation of Changsha (Grant No. kq2208252). The third author was supported by the Excellent Youth Foundation of Education Bureau of Hunan Province (Grant No. 21B0301) and the Natural Science Foundation of Hunan Province (Grant No. 2022JJ40461). And the last author was supported by the Natural Science Foundation of China (Grant No.12101615) and the Natural Science Foundation of Hunan Province (Grant No. 2022JJ40567).





\section*{References}
\bibliographystyle{elsarticle-num}

\end{document}